\documentclass[reqno,12pt]{amsart}
\usepackage{amsmath,amsthm,amscd,amssymb}
\newcommand{\bbC}{{\mathbb{C}}}
\newcommand{\bbD}{{\mathbb{D}}}
\newcommand{\bbR}{{\mathbb{R}}}
\newcommand{\bbZ}{{\mathbb{Z}}}

\newcommand{\lb}{\label}
\newcommand{\f}{\frac}

\newcommand{\ti}{\tilde  }

\newcommand{\loc}{\text{\rm{loc}}}

\newcommand{\ac}{\text{\rm{ac}}}
\newcommand{\s}{\text{\rm{s}}}
\newcommand{\sing}{\text{\rm{sc}}}

\newcommand{\supp}{\text{\rm{supp}}}

\newcommand{\bi}{\bibitem}

\newcommand{\beq}{\begin{equation}}
\newcommand{\eeq}{\end{equation}}
\newcommand{\ba}{\begin{align}}
\newcommand{\ea}{\end{align}}
\newcommand{\veps}{\varepsilon}
\newcounter{smalllist}
\newenvironment{SL}{\begin{list}{{\rm\roman{smalllist})}}{%
\setlength{\topsep}{0mm}\setlength{\parsep}{0mm}\setlength{\itemsep}{0mm}%
\setlength{\labelwidth}{2em}\setlength{\leftmargin}{2em}\usecounter{smalllist}%
}}{\end{list}}

\DeclareMathOperator{\Ima}{Im}

\DeclareMathOperator*{\wlim}{w-lim}

\allowdisplaybreaks
\numberwithin{equation}{section}
\newtheorem{theorem}{Theorem}[section]

\theoremstyle{definition}

\theoremstyle{remark}

\newcommand{\abs}[1]{\lvert#1\rvert}

\begin{document}

\title[Necessary and Sufficient Conditions]
{Necessary and Sufficient Conditions in the Spectral Theory of
Jacobi Matrices and Schr\"odinger Operators}
\author[D.~Damanik, R.~Killip, and B. Simon]{David Damanik$^{1}$,
Rowan Killip$^{2}$, and Barry Simon$^{3}$}

\thanks{$^1$ Mathematics 253-37, California Institute of Technology,
Pasadena, CA 91125, USA. E-mail: damanik@its.caltech.edu. Supported in
part by NSF grant DMS-0227289}
\thanks{$^2$ Department of Mathematics, University of California, Los Angeles,
Los Angeles, CA 90095. E-mail: killip@math.ucla.edu}
\thanks{$^3$ Mathematics 253-37, California Institute of Technology,
Pasadena, CA 91125, USA. E-mail: bsimon@caltech.edu. Supported in
part by NSF grant DMS-0140592}

\date{September 5, 2003}

\begin{abstract} We announce three results in the theory of Jacobi matrices
and Schr\"odinger operators. First, we give necessary and
sufficient conditions for a measure to be the spectral measure of
a Schr\"odinger operator $-\f{d^2}{dx^2} +V(x)$ on $L^2
(0,\infty)$ with $V\in L^2 (0,\infty)$ and $u(0)=0$ boundary
condition. Second, we give necessary and sufficient conditions on
the Jacobi parameters for the associated orthogonal polynomials to
have Szeg\H{o} asymptotics. Finally, we provide necessary and
sufficient conditions on a measure to be the spectral measure of a
Jacobi matrix with exponential decay at a given rate.
\end{abstract}

\maketitle

%%%%%%%%%%%%%%%%%%%%%%%%%%%%%%%
\section{Introduction} \lb{s1}
%%%%%%%%%%%%%%%%%%%%%%%%%%%%%%%

In this note, we want to describe some new results in the spectral
and inverse spectral theory of half-line Schr\"odinger operators
and Jacobi matrices. Given $V\in L_\loc^1 (0,\infty)$ with a mild
regularity condition at infinity (ensuring limit-point case there,
cf.\ \cite{LS}), one can define a unique selfadjoint operator
which is formally
\begin{equation} \lb{1.1}
H=-\f{d^2}{dx^2} + V(x)
\end{equation}
with $u(0) =0$ boundary condition (see, e.g., \cite{LS}). For any
$z\notin\bbR$, there is a solution $u(x;z)$ of $-u'' +Vu =zu$
which is $L^2$ at infinity and unique up to a constant. The Weyl
$m$-function is then defined by
\begin{equation} \lb{1.2}
m(z) = \f{u'(0;z)}{u(0;z)}
\end{equation}
It obeys $\Ima m (z) >0$ when $\Ima z >0$, which implies that
$\Ima m(E+i\veps)$ has a boundary value as $\veps \downarrow 0$ in
distributional sense:
\begin{equation} \lb{1.3}
d\rho(E) = \wlim_{\veps\downarrow 0}\, \f{1}{\pi}\, \Ima m(E+i\veps)\, dE
\end{equation}
$d\rho$ is called the spectral measure.

In this way, each $V$ gives rise to a spectral measure $d\rho$. In
fact, the correspondence is one-to-one: Gel'fand-Levitan
\cite{GLev} (see also Simon \cite{Sim271}) found an inverse
procedure to go from $d\rho$ to $V$\!.

Similarly, given a Jacobi matrix, $a_n >0$, $b_n\in\bbR$:
\begin{equation} \lb{1.4}
J= \begin{pmatrix}
b_1 & a_1 & 0 & 0 & \cdots \\
a_1 & b_2 & a_2 & 0 & \cdots \\
0 & a_2 & b_3 & a_3 & \cdots \\
\vdots & \vdots & \vdots & \vdots & \ddots
\end{pmatrix}
\end{equation}
on $\ell^2 (\bbZ_+)$, we define $d\mu$ to be the measure
associated to the vector $\delta_1$ by the spectral theorem. That
is,
\begin{equation} \lb{1.5}
m(z) \equiv\langle \delta_1, (J-z)^{-1} \delta_1\rangle =\int
\f{d\mu(E)}{E-z}
\end{equation}

In this setting, the inverse procedure dates back to Jacobi,
Chebychev, Markov, and Stieltjes. It is easy to describe: By
applying Gram-Schmidt to $\{1, E, E^2, \ldots \}$ in $L^2(d \mu)$,
we obtain the orthonormal polynomials $p_n(E)$. These obey the
three-term recursion relation
\begin{equation} \lb{1.7}
E p_n (E) =a_{n+1} p_{n+1}(E) + b_{n+1} p_n(E) + a_n p_{n-1}(E)
\end{equation}
Alternatively, one can obtain $a_n,b_n$ from a continued fraction
expansion of $m$ (\cite{Stie,WallBk}).

The main subject of spectral theory is to find relations between
general properties of the spectral measures $d\rho$ or $d\mu$ and
of the differential/difference equation parameters $V$ and
$a_n,b_n$. Clearly, the gems of the subject are ones that provide
necessary and sufficient conditions, that is, a one-to-one
correspondence between some explicit family of measures and some
explicit set of parameters. In this note, we announce three such
results (one involving asymptotics of orthogonal polynomials
rather than the measures) whose details will appear elsewhere
\cite{KS2,Jost1,Jost2}.

In the context of orthogonal polynomials on the unit circle
\cite{SimOPUC}, Verblunsky's form \cite{V36} of Szeg\H{o}'s
theorem \cite{Sz15,Sz20,Sz21} is such a one-to-one correspondence
between a measure and the recurrence coefficients for its
orthogonal polynomials. Baxter's theorem \cite{Bax,Bax2} and
Ibragimov's theorem \cite{Ib,GoIb} can be viewed as other
examples.

Our work here is related to and motivated by the more recent
result of Killip-Simon \cite{KS}:

\begin{theorem}[\cite{KS}]\lb{T1.1} $J-J_0$ is Hilbert-Schmidt, that is
\begin{equation} \lb{1.a}
\sum_{n=1}^\infty (a_n -1)^2 + b_n^2 <\infty
\end{equation}
if and only if the spectral measure $d\mu$ obeys
\begin{SL}
\item[{\rm{(i)}}] {\rm{(Blumenthal-Weyl)}} $\supp (d\mu) = [-2,2] \cup \{E_j^+\}_{j=1}^{N_+} \cup
\{E_j^-\}_{j=1}^{N_-}$ with $E_1^+ > E_2^+ > \cdots > 2$ and $E_1^- < E_2^- < \cdots < -2$ with
$\lim_{j\to\infty} E_j^\pm =\pm2$ if $N_\pm =\infty$.
\item[{\rm{(ii)}}] {\rm{(Normalization)}} $\mu$ is a probability measure.
\item[{\rm{(iii)}}] {\rm{(Lieb-Thirring Bound)}}
\begin{equation}\lb{1.b}
\sum_{\pm, j} (\abs{E_j^\pm} -2)^{3/2}  <\infty
\end{equation}
\item[{\rm{(iv)}}] {\rm{(Quasi-Szeg\H{o} Condition)}}  Let $d\mu_{\ac}(E)=f(E)\, dE$. Then
\begin{equation} \lb{1.c}
\int_{-2}^2 \log [f(E)] \sqrt{4-E^2}\, dE >-\infty
\end{equation}
\end{SL}
\end{theorem}

Our first result is the analog of this theorem for Schr\"odinger
operators. This is discussed in Section~\ref{s2}.

Our second result concerns Szeg\H{o} asymptotics for orthogonal
polynomials. In 1922, Szeg\H{o} \cite{Sz22} proved that if $d\mu
=f(E)\,dE$ where $f$ is supported on $[-2,2]$ and
\begin{equation} \lb{1.8}
\int \log [f(E)] \f{dE}{\sqrt{4-E^2}} > -\infty
\end{equation}
then
\begin{equation} \lb{1.9}
\lim_{n\to\infty}\, z^n p_n (z + z^{-1})
\end{equation}
exists and is nonzero (and finite) for all $z\in\bbD$. There is
work by Gon\v{c}ar \cite{Gon}, Nevai \cite{Nev79}, and Nikishin
\cite{Nik} that allow point masses outside $[-2,2]$. The following
summarizes more recent results on this subject by
Peherstorfer-Yuditskii \cite{PY}, Killip-Simon \cite{KS}, and
Simon-Zlato\v{s} \cite{SZ}:

\begin{theorem}\lb{T1.2} Suppose $d\mu =f(E)\, dE + d\mu_\s$ with $\supp(d\mu_\sing)
\cup\supp(f)\subset [-2,2]$ and
\begin{equation} \lb{1.10}
\sum_{j,\pm}\, (\abs{E_j^\pm}-2)^{1/2} <\infty
\end{equation}
Then the following are equivalent:
\begin{SL}
\item[{\rm{(i)}}] $\inf (a_1 \dots a_n) >0$
\item[{\rm{(ii)}}] All of the following:
\begin{SL}
\item[{\rm{(a)}}]
\begin{equation} \lb{1.d}
\sum_{n=1}^\infty\, \abs{a_n -1}^2 + \abs{b_n}^2 <\infty
\end{equation}
\item[{\rm{(b)}}] $\lim_{n\to\infty} a_n \dots a_1$ exists and is finite and nonzero.
\item[{\rm{(c)}}] $\lim_{n\to\infty} \sum_{j=1}^n b_j$ exists.
\end{SL}
\item[{\rm{(iii)}}]
\begin{equation} \lb{1.e}
\int_{-2}^2 \log [f(E)] \f{dE}{\sqrt{4-E^2}} >-\infty
\end{equation}
\end{SL}

Moreover, if these hold, then the limit \eqref{1.9} exists and is
finite for all $z\in\bbD$ and is nonzero if $z + z^{-1}
\notin\{E_j^\pm\}$.
\end{theorem}

Because \eqref{1.10} is required a priori here, this result is not a necessary
and sufficient condition with only parameter information on one side and only
spectral on the other. In Section~\ref{s3}, we will discuss a necessary and
sufficient condition for the asymptotics \eqref{1.9} to hold, thereby closing
a chapter that began in 1922.

Finally, in Section~\ref{s4}, we discuss necessary and sufficient
conditions on the measure for the $a$'s and $b$'s to obey
\begin{equation} \lb{1.11}
\limsup\, (\abs{a_n-1} + \abs{b_n})^{1/2n} \leq R^{-1}
\end{equation}
for some $R>1$.  Namely, $d\mu$ must give specified weight to
those eigenvalues $E_j$ with $|E_j| < R+R^{-1}$ and the Jost
function must admit an analytic continuation to the disk $\{z :
|z|<R\}$.

The Jost function is naturally defined in terms of scattering;
however, there is a simple procedure for determining it from the
measure and vice versa.  See \eqref{4.2}.

\bigskip

\noindent\textit{Acknowledgment.} We would like to thank Roman
Romanov for drawing our attention to \cite{East}.

%\bigskip
%%%%%%%%%%%%%%%%%%%%%%%%%%%%%%%%%%%%%%%%%%%%%%%%%%%%%%%%%%%%%%%
\section{Schr\"odinger Operators with $L^2$ Potential} \lb{s2}
%%%%%%%%%%%%%%%%%%%%%%%%%%%%%%%%%%%%%%%%%%%%%%%%%%%%%%%%%%%%%%%

The proofs of the results in this section will appear in \cite{KS2}. Given a
measure $d\rho$ on $\bbR$, define $\ti\sigma$ on $[0,\infty)$ by
\begin{equation} \lb{2.1}
\int_0^\infty g\bigl( \sqrt{E}\,\bigr)\, d\rho(E) = \int_0^\infty g(k)\,
d\ti\sigma(k)
\end{equation}
that is, formally $d\ti\sigma(k) = \chi_{(0,\infty)} (k^2)\, d\rho(k^2)$. For the
Schr\"odinger operator with $V=0$,
\begin{equation} \lb{2.2}
d\rho_0 (E) =\pi^{-1} \chi_{[0,\infty)}(E) \sqrt{E}\, dE
\end{equation}
and
\begin{equation} \lb{2.3}
d\ti\sigma_0 (p) =2\pi^{-1} \chi_{[0,\infty)} (p) p^2 \, dp
\end{equation}

Given $\rho$, define $\ti F$ by
\begin{equation} \lb{2.4}
\ti F(q) = \pi^{-1/2} \int_{p\geq 1} p^{-1} e^{-(q-p)^2} [d\ti\sigma(p) -d\ti\sigma_0(p)]
\end{equation}
Since $d\rho$ obeys
\[
\int \f{d\rho(E)}{1+E^2}<\infty
\]
the integral in \eqref{2.4} is convergent.

Define $M(k)$ by $M(k)=m(k^2)$ with $m$ given by \eqref{1.2}. Here is our main result on
$L^2$ potentials:

\begin{theorem}\lb{T2.1} Let $d\rho$ be the spectral measure associated to a potential,
$V$\!. Then $V \in L^2 ([0,\infty))$ if and only if
\begin{SL}
\item[{\rm{(i)}}] {\rm{(Weyl)}} $\supp (d\rho) = [0,\infty) \cup \{E_j\}_{j=1}^N$ with
$E_1 < E_2 < \cdots <0$ and  $\lim_j  E_j =0$ if $N=\infty$.
\item[{\rm{(ii)}}] {\rm{(Local Solubility)}}
\begin{equation} \lb{2.5}
\int_0^\infty \, \abs{\ti F(q)}^2 \, dq <\infty
\end{equation}
\item[{\rm{(iii)}}] {\rm{(Lieb-Thirring)}}
\begin{equation}\lb{2.6}
\sum_j \, \abs{E_j}^{3/2}  <\infty
\end{equation}
\item[{\rm{(iv)}}] {\rm{(Quasi-Szeg\H{o})}}
\begin{equation} \lb{2.7}
\int \log \biggl[ \f{\abs{M(k;0) + ik}^2}{4k\Ima M(k;0)}\biggr] k^2\, dk <\infty
\end{equation}
\end{SL}
\end{theorem}

{\it Remarks.} 1. While there is a parallelism with Theorem~\ref{T1.1}, there are two
significant differences. First, the innocuous normalization condition is replaced with
\eqref{2.5} and, second, \eqref{2.7} involves $M$ and not just $\mu$.

\smallskip
2. Equation~\eqref{2.5} (assuming \eqref{2.7} holds) is an expression of the fact that
$d\rho$ is the spectral measure of an $L_\loc^2$ potential essentially because
it implies (by \cite{GS}) that the $A$-function of \cite{Sim271} is in $L_\loc^2$.

\smallskip
3. The integrand in \eqref{2.7} is $-\log (1-\abs{R(k)}^2)$ where $R$ is a reflection
coefficient. Weak lower semicontinuity of the negative of the entropy used in \cite{KS}
is replaced by lower semicontinuity of the $L^{2n}$ norm.

\smallskip
4. The key to the proof of Theorem~\ref{T2.1} is a strong version
of the Zaharov-Faddeev \cite{ZF} sum rules; essentially following
\cite{KS,SZ,Sim288}, we provide a step-by-step sum rule for $V\in
L_\loc^2$ and take suitable limits. What is interesting is that we
use a whole-line, not half-line, sum rule.

\smallskip
We note that prior to our work, $V\in L^2\Rightarrow$ \eqref{2.6}
was proved by Gardner et al.\ in \cite{GGKM}.  Bounds of this type
are often called Lieb-Thirring inequalities after their work on
moments of eigenvalues for $V\in L^p(\bbR^d)$; see \cite{LT}.
Deift-Killip, \cite{DeiftK}, proved that $V\in L^2$ implies
$f(E)>0$ for a.e.~$E
>0$. There is related work when $-\frac{d^2}{dx^2}+V\geq 0$ in
Sylvester-Winebrenner \cite{SW} and Denisov \cite{DenJDE}.

%%%%%%%%%%%%%%%%%%%%%%%%%%%%%%%%%%%%%%%%
\section{Szeg\H{o} Asymptotics} \lb{s3}
%%%%%%%%%%%%%%%%%%%%%%%%%%%%%%%%%%%%%%%%

The proofs of the results in this section will appear in
\cite{Jost1}. For the study of Szeg\H{o} asymptotics, it is useful
to map $\bbD=\{z : \abs{z}<1\}$ to $\bbC\backslash [-2,2]$ by
$z\to E=z+{z}^{-1}$. Our main result on this issue uses the
following conditions:
\begin{alignat}{2}
&\sum_{n=1}^\infty \, \abs{a_n-1}^2 + \abs{b_n}^2 <\infty \lb{3.1} \\
&\lim_{N\to\infty}\, \sum_{n=1}^N \log (a_n) \qquad&& \text{exists (and is finite)} \lb{3.2} \\
&\lim_{N\to\infty} \, \sum_{n=1}^N b_n \qquad&& \text{exists (and is finite)} \lb{3.3}
\end{alignat}

\begin{theorem}\lb{T3.1} If for some $\veps >0$,
$z^n p_n(z+z^{-1})$ converges uniformly on compact subsets of
$\{z:0<|z|<\varepsilon\}$ to a non-zero value, then
\eqref{3.1}--\eqref{3.3} hold.

Conversely, if \eqref{3.1}--\eqref{3.3} hold, then $z^n
p_n(z+z^{-1})$ converges uniformly on compact subsets of $\bbD$
and has a non-zero limit for those $z\neq0$ where $z+z^{-1}$ is
not an eigenvalue of J.

% and all $z\in\bbC$ with $0<\abs{z}
%<\veps$, the limit \eqref{1.b} exists and is nonzero. Then
%\eqref{3.1}--\eqref{3.3} hold. Conversely, if
%\eqref{3.1}--\eqref{3.3} holds, then for all $z\in\bbD$, the limit
%\eqref{1.b} exists and is nonzero if $z+z^{-1}$ is not an
%eigenvalue of J.
\end{theorem}

{\it Remarks.} 1. By Theorem~\ref{T1.1}, \eqref{3.1} implies only the quasi-Szeg\H{o}
condition \eqref{1.c} whereas all prior discussions of Szeg\H{o} asymptotics have
assumed the stronger Szeg\H{o} condition \eqref{1.e}. We have examples in \cite{Jost1}
where \eqref{3.1}--\eqref{3.3} hold and $\sum (\abs{E_n^\pm}-2)^{1/2}=\infty$ which, by
\cite{SZ}, implies that \eqref{1.e} fails, so we have examples where Szeg\H{o}
asymptotics hold, although the Szeg\H{o} condition fails.

\smallskip
2. The first step in the proof is to show that for fixed
$z\in\bbD$, Szeg\H{o} asymptotics hold if and only if there is a
solution with Jost asymptotics, that is, for which
\begin{equation} \lb{3.4}
\lim z^{-n} u_n(z)
\end{equation}
exists and is non-zero.

\smallskip
3. We have two constructions of the Jost solution when \eqref{3.1}--\eqref{3.3} hold: one
using the nonlocal step-by-step sum rule of \cite{Sim288} and the other using perturbation
determinants \cite{KS}. In either case, one makes a renormalization: In the first
approach, one renormalizes Blaschke products and Poisson-Fatou representations, and, in
the second case, one uses renormalized determinants for Hilbert-Schmidt operators.

\smallskip
While these are the first results we know for Szeg\H{o}/Jost
asymptotics for Jacobi matrices with only $L^2$ conditions,
Hartman \cite{Hart} and Hartman-Wintner \cite{HW} (see also
Eastham \cite[Ch.~1]{East}) have found Jost asymptotics for
Schr\"odinger operators with $V\in L^2$ with conditionally
convergent integral.

\bigskip
%%%%%%%%%%%%%%%%%%%%%%%%%%%%%%%%%%%%%%%%%%%%%%%%%%%%%%%%%%%
\section{Jacobi Parameters With Exponential Decay} \lb{s4}
%%%%%%%%%%%%%%%%%%%%%%%%%%%%%%%%%%%%%%%%%%%%%%%%%%%%%%%%%%%

The proofs of the results in this section will appear in \cite{Jost2}.

If $m$ is given by \eqref{1.5}, we define $M(z)$ by
\begin{equation} \lb{4.1}
M(z) =-m(z+z^{-1})
\end{equation}
Suppose $M(z)$ is the $M$-function of a Jacobi matrix and that
$M(z)$ has an analytic continuation to a neighborhood of
$\bar\bbD$ with the only poles in $\bar\bbD$ lying in
$\bar\bbD\cap\bbR$ and all such poles are simple. Then we can
define
\begin{equation} \lb{4.2}
u(z) =\prod_{k=1}^N b(z,z_k) \exp\biggl( \int \biggl(
\f{e^{i\theta}+z}{e^{i\theta}-z}\biggr) \log \biggl(
\f{\sin\theta}{\Ima M(e^{i\theta})}\biggr)\,
\f{d\theta}{4\pi}\biggr)
\end{equation}
where $\{z_k\}_{k=1}^\infty$ are the poles of $M$ in $\bbD$. This
agrees with the Jost function from scattering theory (see
\cite{KS}), so we call it by this name.

Given $M$ and the Jost function, $u$, suppose $u$ is analytic in
$\{z : \abs{z}<R\}$ and $z_0$ is a zero of $u$ (pole of $M$) with
$R^{-1} <\abs{z_0} <1$. We say $M$ has a canonical weight at $z_0$
if
\begin{equation} \lb{4.3}
\lim_{\substack{z\to z_0 \\ z\neq z_0}}\, (z-z_0) M(z) = (z_0-z_0^{-1}) [u'(z_0) u(z_0^{-1})]^{-1}
\end{equation}

\begin{theorem} \lb{T4.1} Let $M$ be the $M$-function of a Jacobi matrix, $J$. Then $J-J_0$ is
finite rank if and only if
\begin{SL}
\item[{\rm{(i)}}] $M$ is rational and has only simple poles in $\bar\bbD$ with all such
poles in $\bbR$.
\item[{\rm{(ii)}}] $u$ is a polynomial.
\item[{\rm{(iii)}}] $M$ has canonical weight at each $z\in\bbD$ which is a pole of $M$.
\end{SL}
\end{theorem}

\begin{theorem} \lb{T4.2} Let $M$ be the $M$ function of a Jacobi matrix, $J$. Then the
parameters of $J$ obey
\[
(\abs{a_n-1} + \abs{b_n}) \leq C_\veps (R^{-1} +\veps)^{2n}
\]
for some $R>1$ and all $\veps >0$ if and only if
\begin{SL}
\item[{\rm{(i)}}] $M$ is meromorphic on $\{z : \abs{z}<R\}$ with
only simple poles inside $\bar\bbD$ with all such poles in $\bbR$.
\item[{\rm{(ii)}}] $u$ is analytic in $\{z : \abs{z}<R\}$.
\item[{\rm{(iii)}}] $M$ has canonical weight at each pole of $M$,
$z_0\in\bbD$, with $R^{-1} <\abs{z_0}<1$.
\end{SL}
\end{theorem}

Given $u$ and not $m$, there is a normalization issue, so it is easier to discuss the
perturbation determinant \cite{KS} which obeys
\begin{equation} \lb{4.4}
L(z) = \f{u(z)}{u(0)} = \biggl( \, \prod_{n=1}^\infty a_n\biggr) u(z)
\end{equation}

\begin{theorem}\lb{T4.3} Let $L$ be a polynomial with $L(0)=1$. Then $L$ is a perturbation
determinant of a Jacobi matrix, $J$, with $J-J_0$ finite rank if and only if
\begin{SL}
\item[{\rm{(1)}}] $L(z)$ is real if $z\in\bbR$.
\item[{\rm{(2)}}] The only zeros of $L$ in $\bar\bbD$ lie on $\bar\bbD\cap\bbR$ and are simple.
\end{SL}
In that case, there is a unique $J$ with $J-J_0$ finite rank, so $L$ is its perturbation
determinant.
\end{theorem}

{\it Remarks.} 1. While there is a unique $J$ with $J-J_0$ finite rank, if $L$ has $k$ zeros
in $\bbD$, there is a $k$-parameter family of other $J$'s with $L$ as their perturbation
determinant (each such $J$ has $\abs{a_n-1}+\abs{b_n}$ decaying exponentially, but only one
has $J-J_0$ finite rank).

\smallskip
2. There is an analog of Theorem~\ref{T4.3} if $L$ is only
analytic in $\{z : \abs{z}<R\}$.

\smallskip
3. The proofs of these results depend on analyzing the map $(u,M)\to (u^{(1)}, M^{(1)})$ where
$u^{(1)},M^{(1)}$ are the Jost function and $M$-function for $J^{(1)}$, the Jacobi matrix
with the top row and leftmost column removed. We control $||| u^{(1)}|||$ in terms of $|||u|||$
where
\[
|||u|||^2 = \int\, \abs{u(R_1 e^{i\theta}) -u(0)}^2 \, \f{d\theta}{2\pi}
\]
with $R_1 < R$, and are thereby able to show $|||u^{(n)}|||$ goes to zero exponentially.

\smallskip
While we are aware of no prior work on the direct subject of this section, we note that
Nevai-Totik \cite{NT89} solved the analogous problem for orthogonal polynomials on
the unit circle, that Geronimo \cite{Ger} has related results for Jacobi matrices (but
only under an a priori  hypothesis on $M$), and that there are related results in
the Schr\"odinger operator literature (see, e.g., Chadan-Sabatier \cite{ChaSab}).

\bigskip

%%%%%%%%%%%%%%%%%%%%%%%%%%%%%

\end{document}